\documentclass[journal]{IEEEtran}
\IEEEoverridecommandlockouts
\usepackage{comment}
\usepackage{multirow, array}
\usepackage{cite}
\usepackage{eurosym}
\usepackage[mathscr]{euscript}
\usepackage{amsmath,amssymb,amsfonts}
\usepackage{algorithmic}
\usepackage{graphicx}
\usepackage{epstopdf}
\usepackage{textcomp}
\usepackage{tikz}
\usetikzlibrary{arrows,decorations.pathmorphing,backgrounds,fit,positioning,shapes.symbols,chains, shapes,snakes}

\usepackage{xcolor}
\def\BibTeX{{\rm B\kern-.05em{\sc i\kern-.025em b}\kern-.08em
    T\kern-.1667em\lower.7ex\hbox{E}\kern-.125emX}}

\newenvironment{ldescription}[1]
  {\begin{list}{}%
   {\renewcommand\makelabel[1]{##1\hfill}%
   \settowidth\labelwidth{\makelabel{#1}}%
   \setlength\leftmargin{\labelwidth}
   \addtolength\leftmargin{\labelsep}}}
  {\end{list}}

\allowdisplaybreaks

\begin{document}
\title{An Efficient Robust Approach to the Day-ahead Operation of an Aggregator of Electric Vehicles}
\author{\'A. Porras, R. Fern\'andez-Blanco, J. M. Morales and S. Pineda
\thanks{The authors are with the OASYS research group, University of Malaga, Malaga, Spain. E-mail: alvaroporras19@gmail.com; ricardo.fcarramolino@gmail.com; juanmi82mg@gmail.com (corresponding author); spinedamorente@gmail.com}
\thanks{This project has received funding in part by the Spanish Ministry of Science, Innovation and Universities, the Spanish State Research Agency and 
the European Regional Development Fund through project ENE2017-83775-P; and in part by the European Research Council (ERC) under the European Union's Horizon 2020 research and innovation programme (grant agreement No 755705). The authors thankfully acknowledge the computer resources, technical expertise and assistance provided by the SCBI (Supercomputing and Bioinformatics) center of the University of Malaga.}}
\maketitle
\begin{abstract}
The growing use of electric vehicles (EVs) may hinder their integration into the electricity system as well as their efficient operation due to the intrinsic stochasticity associated with their driving patterns. In this work, we assume a profit-maximizer EV-aggregator who participates in the day-ahead electricity market. The aggregator accounts for the technical aspects of each individual EV and the uncertainty in its driving patterns. We propose a hierarchical optimization approach to represent the decision-making of this aggregator. The upper level models the profit-maximizer aggregator's decisions on the EV-fleet operation, while a series of lower-level problems computes the worst-case EV availability profiles in terms of battery draining and energy exchange with the market. Then, this problem can be equivalently transformed into a mixed-integer linear single-level equivalent given the totally unimodular character of the constraint matrices of the lower-level problems and their convexity. Finally, we thoroughly analyze the benefits of the hierarchical model compared to the results from stochastic and deterministic models.
\end{abstract}
\section{Nomenclature} 
The main notation used throughout the text is stated below for quick reference. Symbols $\widehat{\cdot}$ and $\widetilde{\cdot}$ denote expected and realized values, respectively. Other symbols are defined as required. 
\vspace{-0.2cm}
\subsection{Sets and Indices}
\begin{ldescription}{$xxx$}
\item [$\mathcal{T}$] Set of time periods, indexed by $t$.
\item [$\mathcal{V}$] Set of electric vehicles, indexed by $v$.
\end{ldescription}
%
\subsection{Parameters}
\begin{ldescription}{$xxxxxxx$}
\item [$\overline{C}_{v}/\overline{D}_{v}$] Maximum charging/discharging power of electric vehicle $v$ [kW].
\item [$C_{v}^{E}$] Battery cost of electric vehicle $v$ [$\textup{\euro}$/kWh].
\item [$\overline{E}_{v}/\underline{E}_{v}$] Maximum/Minimum energy stored in the battery of electric vehicle $v$ [kWh].
\item [$K_{v}$] Minimum number of periods that electric vehicle $v$ must be available over the time horizon.
\item [$N_T$] Number of time periods of the optimization horizon.
\item[$P^{G}$] Feeder capacity [kW].
\item [$C^p_{1,2}$] Penalty costs [$\textup{\euro}$/kWh].
\item [$S$] Slope of the linear approximation of the battery life as a function of cycles.
\item [$\overline{\alpha}_{v,t}/\underline{\alpha}_{v,t}$]  Maximum/Minimum value of the availability of electric vehicle $v$ in period $t$. 
\item [$\eta_{v}$] Efficiency of electric vehicle $v$.
\item [$\lambda_t$] Day-ahead electricity price in time period $t$ [$\textup{\euro}$/kWh].
\item [$\mathcal{\widehat{\xi}}_{v}$] Expected daily value of the energy required for transportation of electric vehicle $v$ [kWh].
\end{ldescription}
%
\subsection{Variables}
\begin{ldescription}{$xxxxxxl$}
\item [$c_{v,t}/d_{v,t}$] Charging/Discharging power of electric vehicle $v$ in period $t$ [kW].
\item [$c^{D}_{v,t}$] Degradation cost of the battery corresponding to the charging/discharging cycle of  electric vehicle $v$ in period $t$ [\euro].
\item [$e_{v,t}$] Energy stored in the battery of electric vehicle $v$ in period $t$ [kWh].
\item [$p_{t}$] Power exchanged by the aggregator in period $t$ [kW].
\item [$s_{v,t}$] Slack variable used for the energy balance of the battery of electric vehicle $v$ in period $t$ [kWh].
\item [$\alpha_{v,t}$] Availability of electric vehicle $v$ to charge or discharge, being 1 if available and 0 otherwise.
\item [$\psi^{wc}_{v}$] Total net energy injections into the EV-battery associated with the worst-case scenario of the availability profile of electric vehicle $v$ [kWh].
\item [$\tau_{v,t}$] Energy required for transportation of electric vehicle $v$ in period $t$ [kWh].
\end{ldescription}

\section{Introduction}
\label{sec:introduction}
\IEEEPARstart{C}{URRENTLY}, there is a trend towards decentralization due to an increased presence of distributed energy resources (DERs) in power systems. Some important benefits of DERs are: (i) reduced greenhouse gas emissions; (ii) increased power system flexibility; and (iii) increased reliability, resiliency and power quality \cite{contreras2018technical}. Apart from solar and wind generating units, electric vehicles (EVs) are a novel example of DERs currently growing in electricity networks. For practical reasons, the operation of all EVs in a particular area is usually coordinated by an aggregator. According to the recently published European Union's legal framework \cite{directive_2019_944}, an aggregator is defined as \textit{a market participant engaged in aggregation who is not affiliated to the customer's supplier}. In this line, an EV-aggregator can be simply defined as an entity responsible to manage a high number of EVs \cite{guille2009conceptual, lopes2010integration}. The interested reader is referred to \cite{das2020electric} for a thorough review of the main issues and challenges regarding the current situation of the EV market, operational standards and their grid impact, and charging infrastructure, among others. In the electromobility context, the main task of such an aggregator may consist in deciding purchases/sales of electricity to satisfy driving needs at the minimum cost while facing the uncertainty related to its EVs' driving patterns. Therefore, we foresee two main challenges the aggregators may face when dealing with a large EV-fleet, namely the uncertainty in the driving habits of each individual EV while accounting for an accurate representation of their technical and physical characteristics, and, as a consequence, the search of efficient computational methods when scaling to real-sized fleets.

The technical literature includes a wide variety of decision-making models for EV-aggregators that account for uncertainty in driving patterns. A first group of references proposes stochastic optimization problems in which the uncertain parameters are characterized by a finite set of plausible scenarios. Within this category, some works make use of a scenario-based modeling of the uncertainty at the level of the whole EV-fleet, that is, on an aggregate basis \cite{Vagropoulos2013, vaya2015integration, Vaya2016, sarker2017optimal, Rashidizadeh-Kermani2017, Alipour2017a}. Other models, on the contrary, opt for a more refined modeling of the uncertainty through scenarios that capture the electricity demand and availability of each EV in the fleet \cite{iversen2014optimal, momber2014risk, carrion2015operation, Wu2016, hashemi2019stochastic, Habibifar2020}. Although most of these works determine the strategy that maximizes the expected profit of the aggregator, some consider risk-aversion using the conditional value-at-risk \cite{momber2014risk, Alipour2017a, Wu2016, Habibifar2020} or by way of chance-constraints \cite{Vaya2016}. While the decisions yielded by these approaches are more conservative, a complete protection against the worst-case realizations of the uncertainty is not guaranteed. A second disadvantage of the above stochastic models is, besides, the huge number of scenarios required for an accurate representation of the uncertainty, especially if that representation is to go down to the level of the individual EV. In this vein, reference \cite{Habibifar2020} uses a scenario-reduction procedure to trim down the computational burden of the stochastic model they propose.


Alternatively, decision-making models for EV-aggregators based on robust optimization are, comparatively, computationally inexpensive and deliver decisions that maximize the EV-aggregator's profit under the worst possible realization of the uncertainty. Research works in this line are, however,  much scarcer in the technical literature. One instance is reference \cite{battistelli2012optimal}, which considers uncertainty sets for the aggregate power injection of the EV-fleet in vehicle-to-grid (V2G) operational mode. Similarly, the later work \cite{Yang2016} proposes a robust model where an uncertainty set is used for the aggregated energy demand of the fleet. More recently, Baringo and S\'anchez Amaro \cite{baringo2017stochastic} model the behavior of a fleet of EVs as a virtual battery and robust optimization is prescribed to account for the uncertainty in its aggregated power and energy limits. Notice that these references include uncertainty sets for aggregate values of the stochastic variables. In fact, to the best of our knowledge, the technical literature lacks decision-making problems to obtain robust strategies of an EV-aggregator while modeling the uncertain driving patterns of individual EVs.

This paper addresses the day-ahead operation problem of an EV-aggregator who participates in the day-ahead electricity market. This aggregator aims to maximize its profits while scheduling the charging and discharging of each EV in the fleet. We assume that the EVs are equipped with V2G capabilities so that they can operate in two modes, namely grid-to-vehicle (G2V) when extracting power from the grid, and V2G when injecting power into the grid. The aggregator must take into account the physical and technical limitations related to the distribution network (e.g. feeder capacity) and the EVs (e.g. battery degradation, which may hinder the operation of the EV-fleet). However, the source of complexity for this aggregator is to model the uncertain driving patterns of the EVs. Thus, we propose a hierarchical optimization approach that robustifies the operation of the EV-fleet by using past information about the drivers' habits. 
The main contributions of our work are the following:

\begin{itemize}
    \item We propose a novel hierarchical optimization approach for an EV-aggregator to decide the amount of energy to buy or sell in the day-ahead electricity market to cover the driving needs of the EVs in the fleet. 
    \item Since the availability for charging and discharging of the EVs is uncertain, we assume that the aggregator seeks to: i) reduce the risk of the battery of an EV being depleted while driving and ii) decrease the real-time energy deviations of the EV-fleet with respect to the day-ahead plan. For this purpose, our approach simultaneously immunizes the trading strategy of the aggregator against two scenarios for the availability of each EV that are worst-case in terms of battery draining and energy exchange with the market.    
     \item Despite the hierarchical structure of our optimization model, its peculiar characteristics make it efficient from a computational point of view.
\end{itemize}

The result is a robust EV-aggregator's market participation model that is: i) \emph{simple}, in the sense that the uncertain availability of each vehicle is expressed by a few intuitive parameters; ii) \emph{effective}, because it actually hedges the EV-aggregator's trading plan against worst-case EVs' availability patterns; and iii) \emph{efficient}, as it scales well with the number of EVs in the fleet. Furthermore, we thoroughly compare our approach with deterministic and stochastic alternatives on a realistic case study.

This paper builds upon our previous work \cite{8848991}, which we have notably improved and extended. In particular, we consider here EVs with V2G capabilities that, as such, constitute valuable assets for the aggregator to sell energy in the day-ahead market. Due to the inherent risk of battery depletion this action may entail, we rely on the use of two different lower levels to hedge against the EVs' uncertain availability, unlike in \cite{8848991}. Besides, the results are enhanced by comparing them with those from a stochastic model, as well as by performing extensive sensitivity analyses to better illustrate the benefits of the proposed  hierarchical model.

The rest of this paper is organized as follows. Section \ref{sec:formulation} describes three formulations for the decision-making problem of the EV-aggregator: a deterministic approach, a stochastic one, and the proposed hierarchical model.  Section \ref{Methodology} outlines the procedure to transform the original hierarchical problem into a single-level equivalent mixed-integer linear problem (MILP), while Section \ref{Comparison methodologies} explains the methodology we use to benchmark our approach. Section \ref{sec:case} discusses simulation results from a realistic case study. Finally, conclusions are duly drawn in Section \ref{sec:conclusion}.

\section{Problem Formulation}\label{sec:formulation}
In this section, we present the proposed formulation for the day-ahead operation of an EV-aggregator. In Subsection \ref{Deterministic}, we describe the formulation of the decision-making problem faced by the EV-aggregator in its deterministic form. Its stochastic counterpart is presented in Subsection \ref{stochastic_approach}. The characterization of each EV uncertainty and the hierarchical approach are put forward in Subsection \ref{Hierarchical Formulation}. For the sake of unit consistency, we consider hourly periods.
%
\subsection{Deterministic Formulation}
\label{Deterministic}
The EV-aggregator aims to maximize its profits  while both (i) scheduling the charging and discharging of each individual EV and (ii) satisfying physical limitations. Thus, the deterministic formulation can be expressed as:
\begin{subequations}
\label{deterministic_formulation}
\begin{align}
&\min_{\Xi^{D}} \hspace{3pt} \sum_{t \in \mathcal{T}} \widehat{\lambda}_{t}  p_{t} + \sum_{t \in \mathcal{T}}\sum_{v \in \mathcal{V}}  \left(c^{D}_{v,t} + C^p_1 s_{v,t} \right) \label{d1}\\
&\text{subject to:}\notag\\
& p_{t} = \sum_{v \in \mathcal{V}} \left(c_{v,t} - d_{v,t}\right) , \quad \forall t \in \mathcal{T} \label{d2}\\
& -P^{G} \leq p_{t} \leq P^{G}, \quad \forall t \in \mathcal{T} \label{d3}\\
& e_{v,t} = e_{v,t-1} +  \eta_{v} c_{v,t}\hspace{2pt} \widehat{\alpha}_{v,t}  - \frac{d_{v,t}}{\eta_{v}} - \widehat{\tau}_{v,t}+ s_{v,t}, \notag\\
& \hspace{1cm} \forall v \in \mathcal{V}, t \in \mathcal{T} \label{d4}\\
& c_{v,t} \leq \overline{C}_{v}, \quad \forall v \in \mathcal{V}, t \in \mathcal{T} \label{d5}\\
& d_{v,t} \leq \overline{D}_{v} \hspace{2pt} \widehat{\alpha}_{v,t}, \quad \forall v \in \mathcal{V}, t \in \mathcal{T} \label{d6}\\
& \underline{E}_{v} \leq e_{v,t}  \leq \overline{E}_{v}, \quad \forall v \in \mathcal{V}, t \in \mathcal{T} \label{d7}\\
& e_{v,N_T} = e_{v,0}, \quad \forall v \in \mathcal{V}  \label{d8}\\
& c^{D}_{v,t} = \Big{|}\frac{S}{100}\Big{|} C^{E}_{v} \left( \hspace{0.1cm} \frac{1}{\eta_{v}}d_{v,t} + \widehat{\tau}_{v,t} \right),\quad \forall v \in \mathcal{V}, t \in \mathcal{T} \label{d9}\\
& c_{v,t},d_{v,t},e_{v,t},c^{D}_{v,t},s_{v,t} \geq 0, \quad \forall v \in \mathcal{V}, t \in \mathcal{T}, \label{d10}
\end{align}
\end{subequations}
\noindent where the set of decision variables is $\Xi^{D}$ = ($p_{t}$, $c_{v,t}$, $d_{v,t}$, $e_{v,t}$, $c^{D}_{v,t}$, $s_{v,t}$) and parameters $\widehat{\alpha}_{v,t}$ and $\widehat{\tau}_{v,t}$ are expected values.

The objective function \eqref{d1} minimizes the aggregator's operational cost (or equivalently maximizes the aggregator's profit) and it comprises three terms: (a) costs due to the energy bought (i.e., $p_t > 0$) and revenues due to the energy sold (i.e., $p_t < 0$) in the day-ahead market, (b) battery degradation costs, and (c) penalty costs associated with the violation of the evolution of the energy stored in the battery by means of the slack variable $s_{v,t}$. Constraints \eqref{d2} set the power balance of the aggregator at time period $t$. Expressions \eqref{d3} impose the feeder capacity limit at time period $t$. Constraints \eqref{d4} model the evolution of the energy stored in the battery of the EV $v$ between consecutive time periods. If $\widehat{\alpha}_{v,t} = 0$, then the energy of the battery does not increase. Notice that the conversion between electric and chemical energy is affected by a charging and discharging efficiency, which are assumed equal in this paper for the sake of simplicity \cite{ortega2014optimal}. Besides, we introduce slack variables $s_{v,t}$ to model the infeasibilities associated with those EVs that cannot comply with their energy demand. Constraints \eqref{d5} and \eqref{d6} impose the maximum rate of power charging and discharging of EV $v$ at time period $t$, in that order. The maximum and minimum bounds on the energy stored in the battery of the EV $v$ at time period $t$ are set in \eqref{d7}. Constraints \eqref{d8} enforce the terminal condition of the energy stored in the battery. Constraints \eqref{d9} represent the degradation cost  for each charging/discharging cycle of the battery of EV $v$ at time period $t$, where the battery cost $C^{E}_{v}$ represents the purchase cost of the battery divided by its useful capacity. This battery degradation cost model is typically used for Li-ion batteries. Note that expression \eqref{d9} precludes the charging power, because, in the derivation of this approximation, it is presumed that for the EV to be discharged, it had to be previously charged (at some point in the past), so that a charging-discharging cycle could be completed. Consequently, the approximation is formulated as a function of the chemical energy retrieved from the battery, which can be computed as the discharging power divided by the efficiency. Further information on how to estimate this cost can be found in \cite{ortega2014optimal}. Expressions \eqref{d10} define the non-negative character of the decision variables.

In the deterministic model \eqref{deterministic_formulation}, the availability of the EV $v$ at time period $t$, i.e., $\widehat{\alpha}_{v,t}$, is an expected value of a 0/1 random variable (being 0 unavailable and 1 available) and, as such, this expectation can take any value between 0 and 1. 
However, as is customary in everyday life, the availability of an EV is unknown in those periods when the arrival or departure time is uncertain. 

Last but not least, binary variables are not required in problem \eqref{deterministic_formulation} to keep from the simultaneous charging and discharging of the EVs' batteries. This is so because electricity prices are assumed positive, the round-trip efficiency of the battery is strictly smaller than one, the battery degradation costs are accounted for in the objective function \eqref{d1}, and the terminal condition \eqref{d8} is enforced.


\subsection{Stochastic Formulation}
\label{stochastic_approach}

Conversely to the deterministic formulation, this approach characterizes the uncertainty on EVs' availability and consumption through a set of scenarios with given probabilities. Electricity prices are assumed known and equal to their expected values. The stochastic version of the aggregator decision-making model is then formulated as:
\begin{subequations}
\label{stochastic_formulation}
\begin{align}
&\min_{\Xi^{S}} \hspace{3pt}  \sum_{t \in \mathcal{T}}  \widehat{\lambda}_{t}  p_{t} + \sum_{\omega \in \Omega}  \sum_{t \in \mathcal{T}} \sum_{v \in \mathcal{V}} \pi_{\omega} \left( c^{D}_{v,t,\omega} + C^p_1 s_{v,t,\omega}   \right)  \label{s1}\\
&\text{subject to:}\notag\\
& \sum_{v \in \mathcal{V}} \left(c_{v,t,\omega} - d_{v,t,\omega}\right) \leq p_{t}, \quad \forall t \in \mathcal{T},\omega \in \Omega \label{s2}\\
&  -P^{G} \leq p_{t} \leq P^{G}, \quad \forall t \in \mathcal{T} \label{s3}\\
&\left(c_{v,t,\omega}, d_{v,t,\omega}, s_{v,t,\omega}, c^{D}_{v,t,\omega}\right) \in \Phi(\alpha_{v,t,\omega}, \tau_{v,t,\omega}), \notag\\
& \hspace{1cm} \quad \forall v \in \mathcal{V},t \in \mathcal{T},\omega \in \Omega, \label{s4}
\end{align}
\end{subequations}
\noindent wherein $\Omega$ is the set of scenarios, indexed by $\omega$ and $\pi_{\omega}$ is the probability of occurrence of each scenario. The set of decision variables is $\Xi^{S}$ = ($p_{t}$,  $c_{v,t,\omega}$, $d_{v,t,\omega}$, $c^{D}_{v,t,\omega}$, $e_{v,t,\omega}$, $s_{v,t,\omega}$). The objective function minimizes now the average total cost over all scenarios. Importantly, the availability and consumption of each EV, that is, $\alpha_{v,t,\omega}$ and $\tau_{v,t,\omega}$, depend on the scenario realization. In this line, equations \eqref{s4} define the feasible set $\Phi(\cdot)$ represented by constraints \eqref{d4}--\eqref{d10} in terms of the realized values of $\alpha_{v,t}$ and $\tau_{v,t}$ under scenario $\omega$. This approach has two main drawbacks. First, it provides day-ahead strategies that are optimal \textit{on average} and, therefore, low outcomes may occur under adverse realizations of the uncertainty. Second, the high number of variables and constraints may render this problem computationally expensive.

\subsection{Hierarchical Formulation}
\label{Hierarchical Formulation}
\begin{figure}[tbp]
\centerline{\includegraphics[scale=0.35]{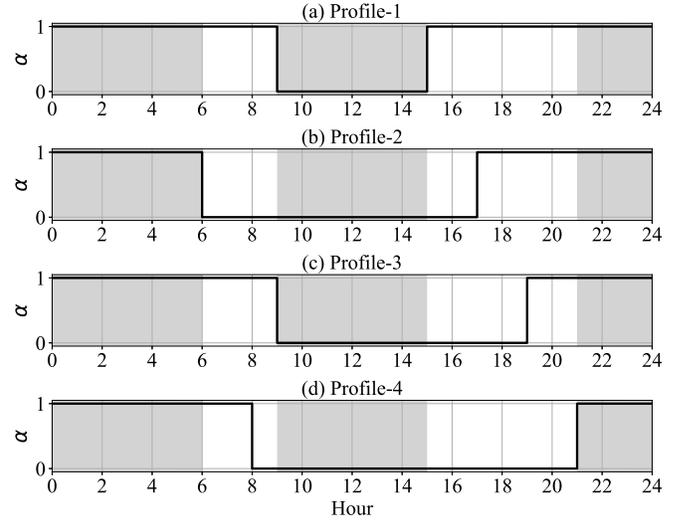}}
\vspace{-0.4cm}
\caption{Four historical profiles of the availability status of an EV. Note that the grey fills represent the time periods when the availability is known.}
\label{fig3}
\end{figure}

As previously mentioned, the time periods when the EVs are available are uncertain in nature. Alternatively to the use of scenarios, we propose to characterize such an uncertainty using the parameter set $\phi_{v}$ = ($K_{v}$, $\overline{\alpha}_{v,t}$, $\underline{\alpha}_{v,t}$), which can be estimated from the corresponding historical availability profiles. Parameter $K_{v}$ stands for the minimum amount of hours within a day that EV $v$ is available. This parameter, however, does not convey \emph{when} exactly the vehicle will be available within the day. This information is conveyed by parameters $\underline{\alpha}_{v,t}$ and  $\overline{\alpha}_{v,t}$. A car for which $\underline{\alpha}_{v,t} \neq \overline{\alpha}_{v,t}$ for many $t$ is a car with a highly uncertain availability pattern. A car with a low $K_{v}$ is a car for which the EV-aggregator has only $K_v$ (potentially unknown) time periods for charging/discharging. By using Fig.~\ref{fig3}, in which we show four historical profiles for an EV, we insist on the meaning of these parameters and clarify how to estimate them.



For instance, let us assume that an EV-aggregator has to determine the optimal strategy for the following day, which happens to be Thursday. Then, a simple, but sensible way to set these uncertain parameters consists in using information of the four previous Thursdays as follows. Let $\kappa_{v,d}$ be the number of hours of day $d$ in which vehicle $v$ is available. Then, for the following day $\check{d}$, parameter $K_v$ can be computed as:
\begin{equation}
K_v = \left\lfloor \frac{1}{4} \left( \kappa_{v,\check{d}-7} + \kappa_{v,\check{d}-14} + \kappa_{v,\check{d}-21} + \kappa_{v,\check{d}-28} \right) \right\rfloor,
\end{equation} 

\noindent where $\lfloor \cdot \rfloor$ is the floor integer function. Similarly, let $\alpha_{v,t}$ be the availability of vehicle $v$ in hour $t$. For each time period $\check{t}$ of the following day, $\overline{\alpha}_{v,\check{t}}$ and $\underline{\alpha}_{v,\check{t}}$ can be determined as:
\begin{align}
& \overline{\alpha}_{v,\check{t}} = 1 - (1-\alpha_{v,\check{t}-168}) \cdot (1-\alpha_{v,\check{t}-336 }) \cdot \notag\\
& \hspace{3.5cm}(1-\alpha_{v,\check{t}-504 }) \cdot (1-\alpha_{v,\check{t}-672 }) \\
& \underline{\alpha}_{v,\check{t}} = \alpha_{v,\check{t}-168} \cdot \alpha_{v,\check{t}-336 } \cdot \alpha_{v,\check{t}-504 } \cdot \alpha_{v,\check{t}-672 }
\end{align}

This is illustrated using Fig. \ref{fig3}, in which we show the availability profiles of a given EV for the previous four Thursdays. We can observe that the EV has been available 18, 13, 14, and 11 periods and, therefore, $K_{v} = 12$. Notice that, in the four profiles, the EV is always available for 00:00-6:00 and 21:00-24:00 and then we can set $\underline{\alpha}_{v,t} = \overline{\alpha}_{v,t} = 1$ for these periods. Similarly, the EV is always unavailable for 9:00-15:00 and then $\underline{\alpha}_{v,t} = \overline{\alpha}_{v,t} = 0$ for these time periods. All these known time periods are indicated in grey fill. In the remaining periods, the EV availability varies depending on the day and, therefore, we set $\underline{\alpha}_{v,t} = 0$ and $\overline{\alpha}_{v,t} = 1$ for 6:00-9:00 and 15:00-21:00. Given the low number of time periods in which the availability is uncertain, this EV can be characterized as a highly predictable EV.

Especially helpful is the fact that our approach encodes the uncertainty in the availability of each individual EV by way of a few intuitive parameters, namely, $K_{v}$, $\overline{\alpha}_{v,t}$, and $\underline{\alpha}_{v,t}$. Actually, in practice, some information on these parameters could be directly provided by the EV users themselves through a home energy management system. For instance, EV users may indicate the time periods of the following day for which their vehicles will be available for charging/discharging with a high level of confidence.



\begin{figure}[t]
\centering
\begin{tikzpicture}
[node distance = 3.4cm, auto,font=\footnotesize,
every node/.style={node distance=2.4cm},
force/.style={rectangle, draw, fill=black!1, inner sep=12pt, text width=6cm, text badly centered, minimum height=1cm, minimum width=8.9cm, font=\bfseries\footnotesize\sffamily},
force2/.style={rectangle, draw, fill=black!1, inner sep=12pt, text width=3.2cm, text badly centered, minimum height=1cm, font=\bfseries\footnotesize\sffamily}]
\node [force] (UL) {Upper-level problem \eqref{h1}--\eqref{h7} \\  (Profit-maximizer aggregator's decisions)};
\node [force2, below of=UL, left of=UL] (LL1) {Lower level \eqref{h8} \\  (Worst-case in terms of battery draining)};
\node [force2, below of=UL, right of=UL] (LL2) {Lower level \eqref{h9} \\  (Worst-case in terms of energy exchange with the market)};
\draw [->,thick] (-2.7,-0.7) -- (-2.7,-1.53) node [midway, left] {$c_{v,t}, d_{v,t}$}  ;
\draw [->,thick] (-1.7,-1.53) -- (-1.7,-0.7) node [midway, right] {$\psi^{wc}_{v}$} ;
\draw [->,thick] (1.7,-0.7) -- (1.7,-1.39) node [midway, left] {$c_{v,t}, d_{v,t}$}  ;
\draw [->,thick] (2.7,-1.39) -- (2.7,-0.7) node [midway, right] {$\alpha_{v, t}$} ;
\end{tikzpicture}

\caption{A conceptual diagram of the interfaces between the upper- and lower-level problems.}
\label{fig:sketch}
\end{figure}
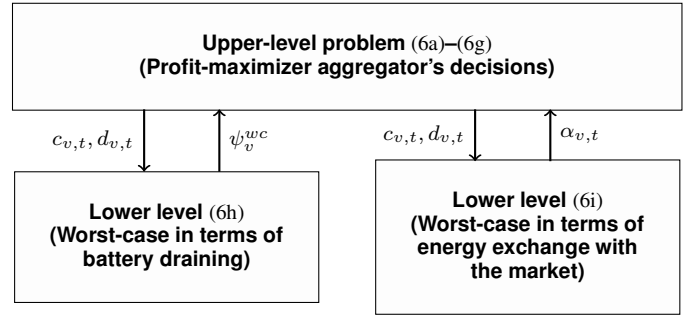

We describe next the proposed approach, which is built upon the deterministic formulation \eqref{deterministic_formulation}. This approach explicitly accounts for the fact that the availability of each EV and, consequently, the energy it requires for transportation are both uncertain and characterized by the parameter set $\phi_{v}$. To do so, the expected values  $\widehat{\tau}_{v,t}$ and $\widehat{\alpha}_{v,t}$ are replaced by decision variables $\tau_{v,t}$ and $\alpha_{v,t}$ in the proposed formulation. The hierarchical optimization model we propose is depicted in Fig.~\ref{fig:sketch}. This figure conceptually represents the interactions among the electricity players, namely the EV-aggregator in the upper level and each of the EVs in the lower level. In addition, the capacity of the feeder may be provided by the community or neighbourhood of the EV-fleet. Specifically, the upper level models the profit-maximizer aggregator's decisions on the EV-fleet operation. The individual charging and discharging plan of the EVs, i.e. $c_{v,t}$ and $d_{v,t}$, are passed on to a series of lower-level problems (two per EV) to compute the EVs availability profiles that are worst-case in terms of battery draining and energy exchange with the market. The outcome of the lower-level problems is the availability profile $\alpha_{v,t}$ and the corresponding total net energy injections into the EV-battery $\psi^{wc}_v$ that results in a robustified operation of the EVs. Both $\alpha_{v,t}$ and $\psi^{wc}_v$ are fed back into the upper-level problem. This approach takes inspiration from the robust contingency-constrained unit commitment proposed in \cite{street2011contingency} and is formulated as the following single-leader-multi-follower problem:
\begin{subequations}
\label{hierarchical_formulation}
\begin{align}
& \hspace{-0.05cm} \min_{\Xi^{R}} \sum_{t \in \mathcal{T}} \widehat{\lambda}_{t}  p_{t} + \sum_{t \in \mathcal{T}}\sum_{v \in \mathcal{V}}  \left(c^{D}_{v,t} + C^p_1 s_{v,t} \right)  \label{h1}\\
& \text{subject to:} \notag\\
& \text{Constraints \eqref{d2}--\eqref{d3}} \label{h2}\\
& (c_{v,t}, d_{v,t}, s_{v,t}, c^{D}_{v,t}) \in \Phi(\alpha_{v,t}, \tau_{v,t}), \quad \forall v \in \mathcal{V}, t \in \mathcal{T} \label{h3}\\
& \sum_{t \in \mathcal{T}} \tau_{v,t} = \widehat{\xi}_{v}, \quad \forall v \in \mathcal{V}  \label{h4}\\
& \tau_{v,t} \leq \left( \overline{E}_{v} - \underline{E}_{v} \right) \left(1 - \alpha_{v,t} \right), \quad \forall v \in \mathcal{V},t \in \mathcal{T}  \label{h5}\\
&\tau_{v,t} \geq 0, \quad \forall v \in \mathcal{V},t \in \mathcal{T}  \label{h6}\\
& \psi^{wc}_{v} \geq  \widehat{\xi}_{v} , \quad \forall v \in \mathcal{V}  \label{h7}\\
& \psi^{wc}_{v} \in \Lambda_v (c_{v,t}, d_{v,t}), \quad \forall v \in \mathcal{V}  \label{h8}\\
& \alpha_{v,t} \in \Upsilon_v (c_{v,t}, d_{v,t}), \quad \forall v \in \mathcal{V}.  \label{h9}
\end{align}
\end{subequations}

The set of decision variables is $\Xi^{R}$ = ($p_{t}$, $c_{v,t}$, $d_{v,t}$, $e_{v,t}$, $\tau_{v,t}$, $c^{D}_{v,t}$, $s_{v,t}$, $\alpha_{v,t}$). Unlike in the deterministic formulation, both $\tau_{v,t}$ and $\alpha_{v,t}$ become now decision variables.

The objective function \eqref{h1} and constraints \eqref{h2} are identical to \eqref{d1}--\eqref{d3}. Expressions \eqref{h3} define the feasible set $\Phi(\cdot)$ represented by constraints \eqref{d4}--\eqref{d10} in terms of the decision variables $\alpha_{v,t}$ and $\tau_{v,t}$. Constraints \eqref{h4} impose that the energy required for transportation throughout the optimization horizon must be equal to the expected daily value demand of EV $v$. The EV energy demand per period, i.e. $\tau_{v,t}$, is now assumed to be unknown for the EV $v$ since it depends on the EV availability, which becomes a decision variable in the proposed formulation. Constraints \eqref{h5} enforce that the expected daily value demand must be distributed among the time periods in which the EV is in a motion status (i.e. $\alpha_{v,t}=0$), being $\left( \overline{E}_{v} - \underline{E}_{v} \right)$ the maximum usable energy capacity of EV-battery $v$ in an hour. Constraints \eqref{h6} define the non-negative character of the decision variables $\tau_{v,t}$. Constraints \eqref{h7} set that the worst case of the total net energy injections into the EV-battery, corresponding to the worst case of availability, must be greater than or equal to the expected daily value of the energy required for transportation. 

As expressed in \eqref{h8}--\eqref{h9}, variables $\psi^{wc}_{v}$ and $\alpha_{v,t}$ are the outcome of two lower-level optimization problems characterizing the availability of the EV $v$, which depend on the upper-level variables $c_{v,t}$ and $d_{v,t}$ given by the aggregator. 
The set $\Lambda_{v}$ in \eqref{h8} includes the availability profiles of the EV $v$ that most severely jeopardize the fulfilling of its expected daily energy demand. The set $\Upsilon_{v}$ in \eqref{h9} comprises, in contrast, those availability profiles that most reduce the interaction of the EV-fleet with the market, thus diminishing the opportunities for making profit. Importantly, the aggregator can plan for a charge greater than the expected daily value of the energy demand per EV. This can be deemed a preventive measure to protect the operation of its EV-fleet against the worst-case availability profiles, thus resulting in a robustified operation of the EV-fleet. Next, we describe the mathematical characterization of the sets $\Lambda_{v}$ and $\Upsilon_{v}$ for each EV $v$. 

\subsubsection{Lower-level Problems Determining the Sets $\Lambda_{v}$}
As mentioned above, the set $\Lambda_{v}$ consists of all those availability profiles that are \emph{worst-case} for the fulfilling of its expected daily energy demand. In other words, these availability profiles are those for which the net energy injection into the EV-battery throughout the optimization horizon is taken to its minimum value. Given the uncertainty set that characterizes the availability profile of the EV in question, these worst-case availability profiles can be computed as the solution to the optimization problem \eqref{lower_level_A} for each EV $v$.
\begin{subequations}
\label{lower_level_A}
\begin{align}
& \psi^{wc}_{v} = \min_{\alpha^{\prime}_{v,t}} \hspace{1pt}  \sum_{t \in \mathcal{T}}  \hspace{1pt}  \alpha^{\prime}_{v,t}\left(\eta_{v} c_{v,t} - \frac{1}{\eta_{v}} d_{v,t}\right) \label{ll1-1}\\
&
\hspace{1.1cm} \text{subject to:} \notag\\
& \hspace{1.1cm}\sum_{t \in \mathcal{T}} \alpha^{\prime}_{v,t} \geq K_{v} : (\zeta_{v}^{\prime})  \label{ll1-2}\\
& \hspace{1.1cm} \underline{\alpha}_{v,t} \leq \alpha^{\prime}_{v,t} \leq  \overline{\alpha}_{v,t} : (\underline{\beta}_{v,t}^{\prime}, \overline{\beta}_{v,t}^{\prime}),\quad \forall t \in \mathcal{T} \label{ll1-3}\\
& \hspace{1.1cm} \alpha^{\prime}_{v,t} \in \{0, 1\},\quad \forall t \in \mathcal{T}, \label{ll1-4}
\end{align}
\end{subequations}
\noindent where $\alpha^{\prime}_{v,t}$ are the decision variables. Dual variables for the relaxed linear version of \eqref{lower_level_A} are shown in parentheses after a colon next to the corresponding constraint. Problem \eqref{lower_level_A} is parametrized in terms of the upper-level variables $c_{v,t}$ and $d_{v,t}$.

The objective function \eqref{ll1-1} aims to minimize the chemical energy stored in the battery over the day computed through the injection/retrieval of electric energy and the charging/discharging efficiency. Constraint \eqref{ll1-2} sets the minimum number of time periods ($K_{v}$) that the EV $v$ must be available throughout the time horizon. Expressions \eqref{ll1-3} set the minimum and maximum availability status and serve us to easily enforce the availability or unavailability of the EV $v$ in time period $t$ by imposing $\underline{\alpha}_{v,t} = \overline{\alpha}_{v,t} = 1$ or $\underline{\alpha}_{v,t} = \overline{\alpha}_{v,t} = 0$, respectively. Constraints \eqref{ll1-4} set the binary character of variables $\alpha^{\prime}_{v,t}$. Naturally, the parameter set $\phi_{v}$ $=$ ($K_v$, $\overline{\alpha}_{v,t}$, $\underline{\alpha}_{v,t}$) should be estimated \textit{a priori} from historical data, as discussed above. 

\begin{figure}[tbp]
\centerline{\includegraphics[scale=0.35]{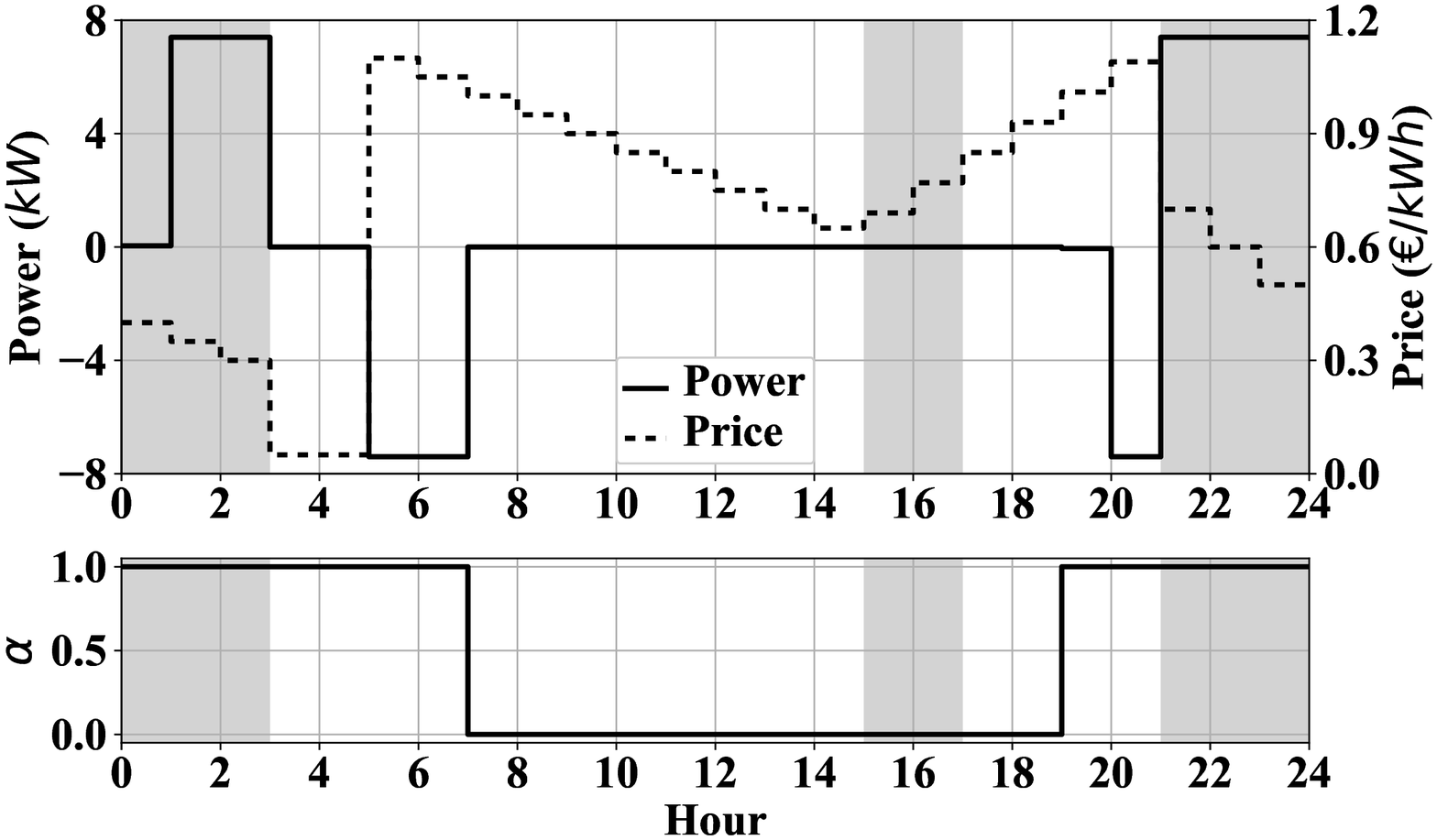}}
\vspace{-0.4cm}
\caption{Schedule for an EV when using the hierarchical optimization model \eqref{h1}--\eqref{h8}, i.e. with the lower level $\Lambda_{v}$ only. Note that grey fills are the time periods in which the availability is known.}
\label{fig5}
\end{figure}

Problem \eqref{lower_level_A} essentially aims to drain the battery of the EV as much as possible. Note that, variables $c_{v,t}$ and $d_{v,t}$ cannot take simultaneously non-zero values for the EV $v$ at time period $t$. Then, the availability is pushed to be 0 when the coefficient of the objective function \eqref{ll1-1}, i.e. $\eta_{v} c_{v,t} - \frac{1}{\eta_{v}} d_{v,t}$, is positive, and 1 when it is negative. However, the worst case of availability profiles that problem \eqref{lower_level_A} determines for each EV does not deter the aggregator from planning large charging-discharging cycles of the EV-batteries in order to perform market arbitrage, as illustrated in Fig.~\ref{fig5}. 

The upper plot of Fig.~\ref{fig5} shows the charging/discharging plan for a certain EV, as given by the hierarchical model \eqref{h1}--\eqref{h8} once transformed into an MILP according to the methodology outlined in \cite{8848991}. The lower plot depicts the availability profile that leads to the maximum draining of the EV-battery. Such a profile does not keep the aggregator from implementing an aggressive market strategy essentially based on arbitrage, which may potentially result in an unacceptably risky operation of the EV-fleet. In order to obtain more conservative strategies, the considered availability profile is obtained instead through set $\Upsilon_{v}$.

\subsubsection{Lower-level Problems Determining the Sets $\Upsilon_{v}$}
The lower levels presented in the previous subsection do not prevent the aggregator from planning the discharge of an EV (in the hope of benefiting from market arbitrage) in time periods where the availability of that EV is uncertain. To cope with this issue, we force the scheduling plan of the aggregator to be also robust against a second type of worst-case availability profiles, namely those whereby the interaction of the EV with the grid, and thus, with the market is minimized. Consequently, the consideration of these availability profiles encourages the aggregator to play a less aggressive but safer market strategy. The set $\Upsilon_{v}$ for each EV $v$ is mathematically given by the following optimization problem:
\begin{subequations}
\label{lower_level_B}
\begin{align}
&  \min_{\alpha_{v,t}} \hspace{1pt} \sum_{t \in \mathcal{T}} \hspace{1pt}  \alpha_{v,t}\left(\eta_{v} c_{v,t} + \frac{1}{\eta_{v}} d_{v,t}\right) \label{ll2-1}\\
&
\text{subject to:} \notag\\
& \sum_{t \in \mathcal{T}} \alpha_{v,t} \geq  K_{v} : (\zeta_{v})  \label{ll2-2}\\
& \underline{\alpha}_{v,t} \leq \alpha_{v,t} \leq  \overline{\alpha}_{v,t}: (\underline{\beta}_{v,t}, \overline{\beta}_{v,t}),\quad \forall t \in \mathcal{T} \label{ll2-3}\\
& \alpha_{v,t} \in \{0, 1\},\quad \forall t \in \mathcal{T}, \label{ll2-4}
\end{align}
\end{subequations}
\noindent where $\alpha_{v,t}$ are the decision variables and dual variables for the relaxed linear version of problem \eqref{lower_level_B} are shown in parentheses after a colon next to the corresponding constraint. Problem \eqref{lower_level_B} is structurally almost identical to optimization problem \eqref{lower_level_A}, except for the objective function, which, in this case, minimizes the interaction of the EV with the market. In other words, problem \eqref{lower_level_B} determines the availability profile that minimizes the planned exchange of power between the EV and the main power grid through the feeder. Note that lower level \eqref{lower_level_B} intends to push the availability statuses of the EV to be zero when the aggregator schedules either charging or discharging, while also satisfying the expected energy required for transportation \eqref{h7}. 

\begin{figure}[tbp]
\centerline{\includegraphics[scale=0.35]{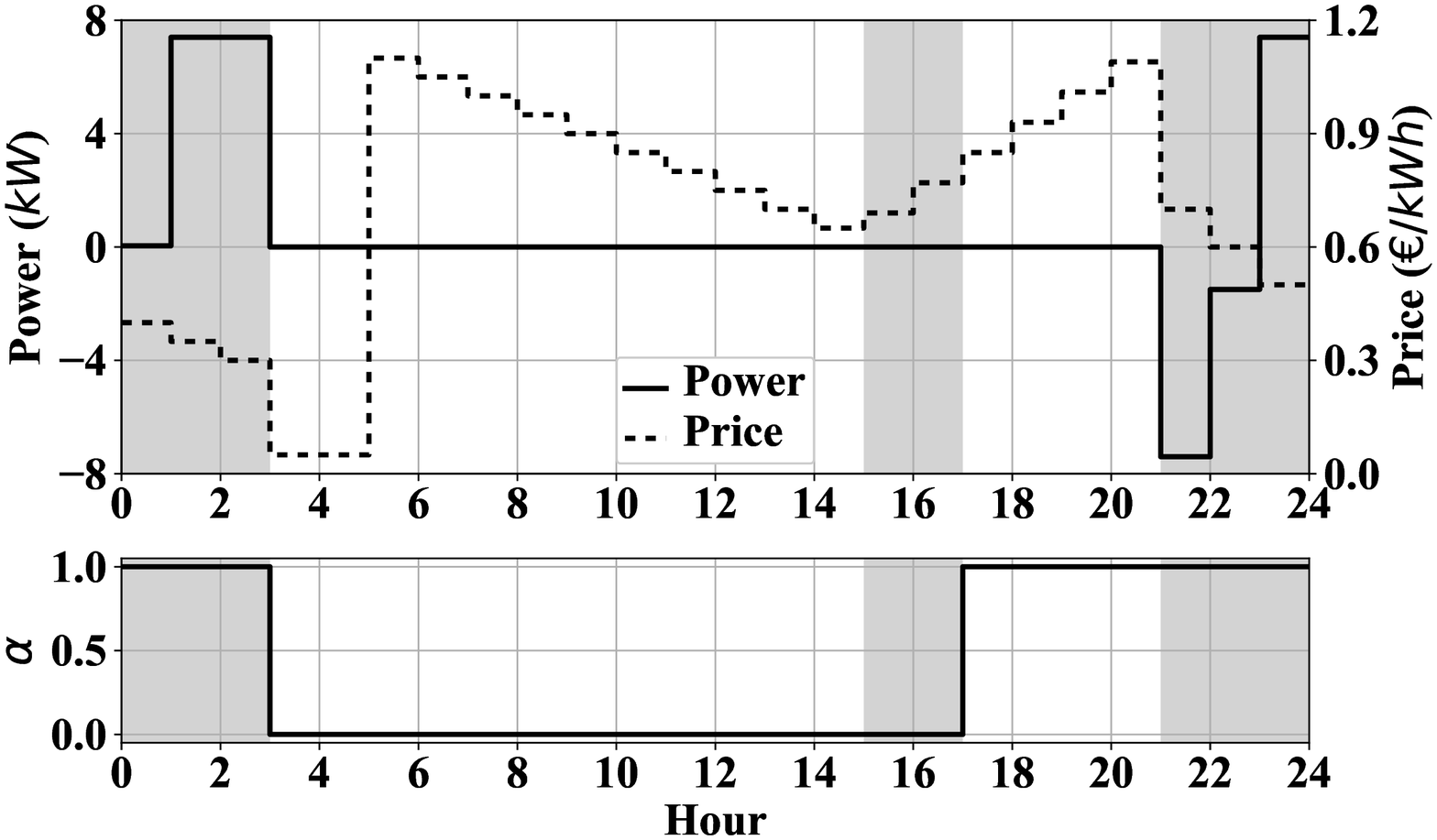}}
\vspace{-0.4cm}
\caption{Schedule for an EV when using the hierarchical optimization model \eqref{hierarchical_formulation}, i.e. with both lower levels, $\Lambda_{v}$ and $\Upsilon_{v}$. Note that grey fills are the time periods in which the availability is known.}
\label{fig6}
\end{figure}

To better illustrate the impact of problem \eqref{lower_level_B}, Fig. \ref{fig6} represents the solution of the EV used to obtain the results given in Fig. \ref{fig5}. To do that, we solve the hierarchical model \eqref{hierarchical_formulation}, i.e. with both lower levels \eqref{lower_level_A} and \eqref{lower_level_B}, once transformed into an MILP according to the methodology described in Section \ref{Methodology}. As can be seen, the aggregator schedules a charging and a discharging plan for the EV in the time periods when the availability is known, i.e. periods 1--3, 22--24. The aggregator also obtains profits, but they are lower than when using only the lower level \eqref{lower_level_A} for the EV in question. Therefore, the use of both lower-level problems  leads to a more conservative decision than when using just the lower level \eqref{lower_level_A}.

\section{Methodology}
\label{Methodology}
The hierarchical problem \eqref{hierarchical_formulation} can be transformed into a mixed-integer nonlinear single-level equivalent as follows: 
\begin{enumerate}
    \item We first verify that $\sum_{t\in\mathcal{T}} \overline{\alpha}_{v,t} \geq K_v$, so that problems \eqref{lower_level_A} and \eqref{lower_level_B} are guaranteed to have, at least, the feasible solution $\alpha_{v,t}=\overline{\alpha}_{v,t},  \forall t\in\mathcal{T}$.
    \item Lower-level problems \eqref{lower_level_A} and \eqref{lower_level_B} are non-convex because variables $\alpha_{v,t}^{\prime}$ and $\alpha_{v,t}$ are binary. We can relax problems \eqref{lower_level_A} and \eqref{lower_level_B} by removing constraints \eqref{ll1-4} and \eqref{ll2-4}. By doing so, $\alpha_{v,t}^{\prime}$ and $\alpha_{v,t}$ can now take any value between 0 and 1 for uncertain periods for which $\underline{\alpha}_{v,t}=0$ and $\overline{\alpha}_{v,t}=1$.

Optimization problems \eqref{lower_level_A} and \eqref{lower_level_B} satisfy that: (i) their corresponding constraint matrices are totally unimodular (TU), as we show at the end of this section, and (ii) the parameters $K_v$, $\overline{\alpha}_{v,t}$, and $\underline{\alpha}_{v,t}$ characterizing the EVs are integer. Under these conditions, we can guarantee that there exists an optimal solution to the relaxed lower-level problems for which the variables $\alpha_{v,t}^{\prime}$ and $\alpha_{v,t}$ take integer values. Therefore, the lower levels \eqref{lower_level_A} and \eqref{lower_level_B} can be replaced by their relaxed formulations.
    \item Under the assumption of convexity of the relaxed lower-level problems \eqref{lower_level_A} and \eqref{lower_level_B}, results from duality theory of linear programming \cite{luenberger} can be applied to transform the original hierarchical program into a non-linear single-level equivalent as follows:
    \begin{enumerate}
        \item $\psi^{wc}_{v}$ can be replaced by the dual objective function of the lower-level optimization problem \eqref{ll1-1}--\eqref{ll1-3} for each EV $v$.
        \item The lower level \eqref{lower_level_A} can be then replaced by the dual feasibility constraints associated with \eqref{ll1-1}--\eqref{ll1-3}. 
        \item The lower level \eqref{lower_level_B} can be replaced by its primal and dual feasibility constraints as well as the equality corresponding to the strong duality condition.
    \end{enumerate}
    \item In order to avoid multiplicity of solutions because of the continuous character of variables $\alpha_{v,t}$, we enforce back its binary character in the single-level equivalent. 
\end{enumerate}
The resulting mixed-integer nonlinear single-level equivalent can be written as follows: 
\begin{subequations}
\label{nonlinear_singlelevel_equivalent}
\begin{align}
&\min_{\Xi^{R}} \hspace{3pt} \sum_{t \in \mathcal{T}} \widehat{\lambda}_{t}  p_{t} + \sum_{t \in \mathcal{T}} \sum_{v \in \mathcal{V}}  \left(c^{D}_{v,t} + C^p_1 s_{v,t}\right)  \label{sl1}\\
&\text{subject to:}\notag\\
& p_{t} = \sum_{v \in \mathcal{V}} \left(c_{v,t} - d_{v,t}\right) , \quad \forall t \in \mathcal{T} \label{sl2}\\
& -P^{G} \leq p_{t} \leq P^{G}, \quad \forall t \in \mathcal{T} \label{sl3}\\
& e_{v,t} = e_{v,t-1} + \eta_{v} c_{v,t} \alpha_{v,t}  - \frac{d_{v,t}}{\eta_{v}}  - \tau_{v,t} + s_{v,t} , \notag\\
& \hspace{2.5cm} \forall v \in \mathcal{V}, t \in \mathcal{T} \label{sl4}\\
& c_{v,t} \leq \overline{C}_{v} , \quad \forall v \in \mathcal{V}, t \in \mathcal{T} \label{sl5}\\
& d_{v,t} \leq \overline{D}_{v} \alpha_{v,t}, \quad \forall v \in \mathcal{V}, t \in \mathcal{T} \label{sl6}\\
& \underline{E}_{v} \leq e_{v,t}  \leq \overline{E}_{v}, \quad \forall v \in \mathcal{V}, t \in \mathcal{T} \label{sl7}\\
& e_{v,N_T} = e_{v,0}, \quad \forall v \in \mathcal{V}  \label{sl8}\\
& c^{D}_{v,t} = \Big{|}\frac{S}{100}\Big{|} C^{E}_{v} \left( \hspace{0.1cm} \frac{1}{\eta_{v}}d_{v,t} + \tau_{v,t} \right),\quad \forall v \in \mathcal{V}, t \in \mathcal{T} \label{sl9}\\
& c_{v,t},d_{v,t},e_{v,t},c^{D}_{v,t},s_{v,t} \geq 0, \quad \forall v \in \mathcal{V}, t \in \mathcal{T} \label{sl10}\\
& \sum_{t \in \mathcal{T}} \tau_{v,t} = \widehat{\xi}_{v}, \quad \forall v \in \mathcal{V}  \label{sl11}\\
& \tau_{v,t} \leq \left( \overline{E}_{v} - \underline{E}_{v} \right) \left(1 - \alpha_{v,t} \right), \quad \forall v \in \mathcal{V},t \in \mathcal{T}  \label{sl12}\\
&\tau_{v,t} \geq 0, \quad \forall v \in \mathcal{V},t \in \mathcal{T}  \label{sl13}\\
&  K_{v} \zeta_{v}^{\prime} + \sum_{t \in \mathcal{T}} \left(\underline{\alpha}_{v,t} \underline{\beta}_{v,t}^{\prime}  + \overline{\alpha}_{v,t} \overline{\beta}_{v,t}^{\prime} \right) \geq  \widehat{\xi}_{v}, \quad \forall v \in \mathcal{V} \label{sl14} \\
& \zeta_{v}^{\prime} + \underline{\beta}_{v,t}^{\prime} + \overline{\beta}_{v,t}^{\prime} =  \eta_{v} c_{v,t} - \frac{1}{\eta_{v}} d_{v,t},  \quad \forall v \in \mathcal{V}, t \in \mathcal{T} \label{sl15}\\
&  \underline{\beta}_{v,t}^{\prime} \geq 0, \overline{\beta}_{v,t}^{\prime} \leq 0, \quad \forall v \in \mathcal{V}, t \in \mathcal{T} \label{sl16}\\
& \zeta_{v}^{\prime} \geq 0, \quad \forall v \in \mathcal{V} \label{sl17}\\
& \sum_{t \in \mathcal{T}} \alpha_{v,t} \geq K_{v}, \quad \forall v \in \mathcal{V}  \label{sl18}\\
& \underline{\alpha}_{v,t} \leq \alpha_{v,t} \leq  \overline{\alpha}_{v,t},\quad \forall v \in \mathcal{V}, t \in \mathcal{T} \label{sl19}\\
& \zeta_{v} + \underline{\beta}_{v,t} + \overline{\beta}_{v,t} =  \eta_{v} c_{v,t} + \frac{1}{\eta_{v}} d_{v,t},  \quad \forall v \in \mathcal{V}, t \in \mathcal{T} \label{sl20}\\
&  \underline{\beta}_{v,t} \geq 0, \overline{\beta}_{v,t} \leq 0, \quad \forall v \in \mathcal{V}, t \in \mathcal{T} \label{sl21}\\
& \zeta_{v} \geq 0, \quad \forall v \in \mathcal{V} \label{sl22}\\
& K_{v} \zeta_{v} + \sum_{t \in \mathcal{T}} \left(\underline{\alpha}_{v,t} \underline{\beta}_{v,t}  + \overline{\alpha}_{v,t} \overline{\beta}_{v,t} \right)  = \notag\\
&\hspace{2cm} \sum_{t \in \mathcal{T}} \alpha_{v,t}\left(\eta_{v} c_{v,t} + \frac{1}{\eta_{v}} d_{v,t}\right), \quad \forall v \in \mathcal{V} \label{sl23}\\
& \alpha_{v,t} \in \{0, 1\}, \quad \forall t \in \mathcal{T}, v \in \mathcal{V}. \label{sl24}
\end{align}
\end{subequations}
The objective function \eqref{sl1} is identical to the upper-level objective function \eqref{h1}. For the sake of completeness, we fully show constraints \eqref{sl2}--\eqref{sl13}, which correspond to the upper-level constraints \eqref{h2}--\eqref{h6}. Expressions \eqref{sl14}--\eqref{sl17} are equivalent to constraints \eqref{h7} and the lower level \eqref{h8}, in which the dual objective function of the lower level \eqref{h8} replaces the variables $\psi^{wc}_{v}$ in constraints \eqref{sl14}. Expressions \eqref{sl15}--\eqref{sl17} are the constraints of the dual of the lower level \eqref{h8}. Constraints \eqref{sl18}--\eqref{sl23} replace the lower level \eqref{h9} with its primal feasibility constraints \eqref{sl18}--\eqref{sl19}, its dual feasibility constraints \eqref{sl20}--\eqref{sl22}, and the strong duality condition \eqref{sl23}. Expression \eqref{sl24} defines the binary character of variables $\alpha_{v,t}$. The single-level equivalent \eqref{nonlinear_singlelevel_equivalent} is characterized as a nonlinear mixed-integer program due to the presence of nonlinear products between bounded continuous variables and binary variables (i.e., $c_{v,t} \alpha_{v,t}$ and $d_{v,t} \alpha_{v,t}$ in constraints \eqref{sl4} and \eqref{sl23}), so they can be easily linearized by using integer algebra results \cite{floudas1995}. Therefore, constraints \eqref{sl4} and \eqref{sl23} can be replaced with the following set of linear constraints:
\begin{subequations}
\begin{align}
& e_{v,t} \hspace{-1.2pt}=\hspace{-1.2pt} e_{v,t-1}\hspace{-1.2pt} + \hspace{-1.2pt}\eta_{v} z^{c}_{v,t}\hspace{-1.2pt} -\hspace{-1.2pt} \frac{ d_{v,t}}{\eta_{v}}\hspace{-1.2pt} -\hspace{-1.2pt} \tau_{v,t} \hspace{-1.2pt}+ \hspace{-1.2pt}s_{v,t}, \forall t\hspace{-1.2pt} \in\hspace{-1.2pt} \mathcal{T}, v \hspace{-1.2pt}\in\hspace{-1.2pt} \mathcal{V} \label{l1}\\
& K_{v} \zeta_{v} + \sum_{t \in \mathcal{T}} \left(\underline{\alpha}_{v,t} \underline{\beta}_{v,t}  + \overline{\alpha}_{v,t} \overline{\beta}_{v,t} \right) =  \notag\\
&\hspace{2.5cm}  \sum_{t \in \mathcal{T}} \left( \eta_{v} z^{c}_{v,t} + \frac{z^{d}_{v,t}}{\eta_{v}} \right), \quad   \forall v \in \mathcal{V} \label{l2}\\
& 0 \leq c_{v,t} - z^{c}_{v,t} \leq \left( 1 - \alpha_{v,t} \right) \overline{C}_v, \quad \forall t \in \mathcal{T}, v \in \mathcal{V} \label{l3}\\
& 0 \leq z^{c}_{v,t} \leq \alpha_{v,t} \overline{C}_v, \quad \forall t \in \mathcal{T}, v \in \mathcal{V}. \label{l4}\\
& 0 \leq d_{v,t} - z^{d}_{v,t} \leq \left( 1 - \alpha_{v,t} \right) \overline{D}_v, \quad \forall t \in \mathcal{T}, v \in \mathcal{V} \label{l5}\\
& 0 \leq z^{d}_{v,t} \leq \alpha_{v,t} \overline{D}_v, \quad \forall t \in \mathcal{T}, v \in \mathcal{V}. \label{l6}
\end{align}
\end{subequations}

The set of decision variables is $\Xi^{R}$ = ($p_{t}$, $c_{v,t}$, $d_{v,t}$, $e_{v,t}$, $\tau_{v,t}$, $c^{D}_{v,t}$, $z^{c}_{v,t}$, $z^{d}_{v,t}$, $s_{v,t}$, $\alpha_{v,t}$, $\overline{\beta}^{\prime}_{v,t}$, $\underline{\beta}^{\prime}_{v,t}$, $\zeta^{\prime}_{v}$, $\overline{\beta}_{v,t}$, $\underline{\beta}_{v,t}$, $\zeta_{v}$).

As previously mentioned, the relaxed lower-level problems \eqref{lower_level_A} and \eqref{lower_level_B} differ from each other only in their objective function. They share, therefore, the same constraint matrix. There are several ways to prove that this matrix is TU. If we first resort to Proposition~2.3 in \cite{wolsey2014integer}, Chapter III.1, the proof boils down to showing that the matrix formed by all the rows of the constraint matrix of the lower-level problems except for those corresponding to bounds on the $\alpha$-variables is TU. This matrix turns out to be the one-row matrix $(1, 1, \ldots, 1)$, whose total unimodularity trivially follows from the fact that $(1, 1, \ldots, 1)$ is a $(0,1,-1)$-matrix with only one single nonzero entry in each column (see Proposition~2.6 in \cite{wolsey2014integer}, Chapter III.1). 

\section{Comparison Methodologies}
\label{Comparison methodologies}
The results from the proposed hierarchical formulation, which is hereinafter referred to as HF, are compared against two benchmarks: (i) Deterministic Formulation (see Subsection \ref{Deterministic}) and (ii) Stochastic Formulation (see Subsection \ref{stochastic_approach}). These methods are denoted as DF and SF, respectively. 

All these three approaches make use of the \emph{same} data: Historical availability profiles for the EVs in the fleet. In this respect, they differ in the type of information they need to retrieve from this data. More particularly, the deterministic model requires expectations of $\alpha_{v,t}$ only, while the stochastic one necessitates a number of scenarios or path-time trajectories for these stochastic processes. For its part, our approach requires the estimation of $K_v$, $\overline{\alpha}_{v,t}$ and $\underline{\alpha}_{v,t}$ from that data. Note that, in terms of estimation effort or complexity, computing expectations only is probably easier than determining $K_v$, $\overline{\alpha}_{v,t}$ and $\underline{\alpha}_{v,t}$ (these parameters can be interpreted as quantiles of appropriate probability distributions), which, in turn, is easier than fabricating a whole set of scenarios.

On a different front, to fairly compare the performance of these methods, we solve a feasibility problem by assuming that the uncertainty is realized and the aggregator must at least sell the energy agreed in the day-ahead market (hereinafter denoted as the \emph{minimum value of the energy sold}) and consume, at most, the energy purchased in such a market. In this setup, we can minimize the magnitude of the deviations from the energy balance of EV-batteries to satisfy their energy demand and the deviations caused by unfulfilling the minimum value of the energy sold, as given by expression \eqref{rt1}. Thus, this feasibility problem can be formulated as:
\begin{subequations}
\label{feasibility_problem}
\begin{align}
&\min_{\Xi^{FP}} \hspace{3pt} \sum_{t \in \mathcal{T}} C^p_2 \hspace{1pt} p^{-}_{t} + \sum_{t \in \mathcal{T}} \sum_{v \in \mathcal{V}}  C^p_1 \hspace{1pt} s_{v,t} \label{rt1}\\
&\text{subject to:}\notag\\
&  \sum_{v \in \mathcal{V}} \left( c_{v,t} - d_{v,t} \right) \leq p_t + \mathbb{I}\{p_t<0\} \cdot p_t^{-}, \quad \forall t \in \mathcal{T} \label{rt2}\\
& \left(c_{v,t}, d_{v,t}, s_{v,t}\right) \in \Phi(\widetilde{\alpha}_{v,t}, \widetilde{\tau}_{v,t}), \quad \forall v \in \mathcal{V},  t \in \mathcal{T} \label{rt10} \\
& p^{-}_{t} \geq 0, \quad \forall  t \in \mathcal{T} \label{rt11},
\end{align}
\end{subequations}
\noindent where the set of decision variables $\Xi^{FP}$ = ($c_{v,t}$, $d_{v,t}$, $e_{v,t}$, $s_{v,t}$, $p^-_{t}$) and $\mathbb{I}\{p_t<0\}$ is the indicator operator, being $1$ if $p_t < 0$, and $0$ otherwise. Constants $C^p_1$ and $C^p_2$ are penalty parameters.

Constraints \eqref{rt2} set the maximum charging of the EV-fleet as the power bought in the day-ahead market and its minimum discharging as the power sold in the day-ahead market. Constraints \eqref{rt10} define the feasibility set $\Phi(\cdot)$, given by constraints \eqref{d4}--\eqref{d10}, in terms of the actual realizations of the availability and consumption, i.e. $\widetilde{\alpha}_{v,t}$ and $\widetilde{\tau}_{v,t}$. Note that there is no need to take the battery degradation costs into account in this feasibility problem. Finally, expressions \eqref{rt11} impose the non-negativity character of variables $p^{-}_t$. In this feasibility problem, slack variables $s_{v,t}$ are the violations of the energy required for transportation of the EV $v$ at time period $t$. In addition, slack variables $p^{-}_{t}$ represent the violation of the power sold at time period $t$. Since the main purpose of the aggregator is to satisfy the EVs' energy demand for transportation, we assume, for the sake of comparison, that the penalty parameter $C^p_1 >> C^p_2$.

\section{Case Study}\label{sec:case}
We analyze here the benefits of the proposed approach through simulation results. In Subsection \ref{Data}, the data for a real-life case study are presented. In Subsection \ref{Results for the Base Case}, we analyze the results for a base case. Next, we provide a sensitivity analysis of the robustness of our model to the 
parameter $K_{v}$ in Subsection \ref{Impact of the Parameter K_v}. Finally, in Subsection \ref{Impact of the Feeder Capacity}, we show scalability results for a fleet of 1000 EVs including a sensitivity analysis on the feeder capacity.

\subsection{Data}
\label{Data}
\begin{figure}[t]
    \centerline{\includegraphics[scale=0.35]{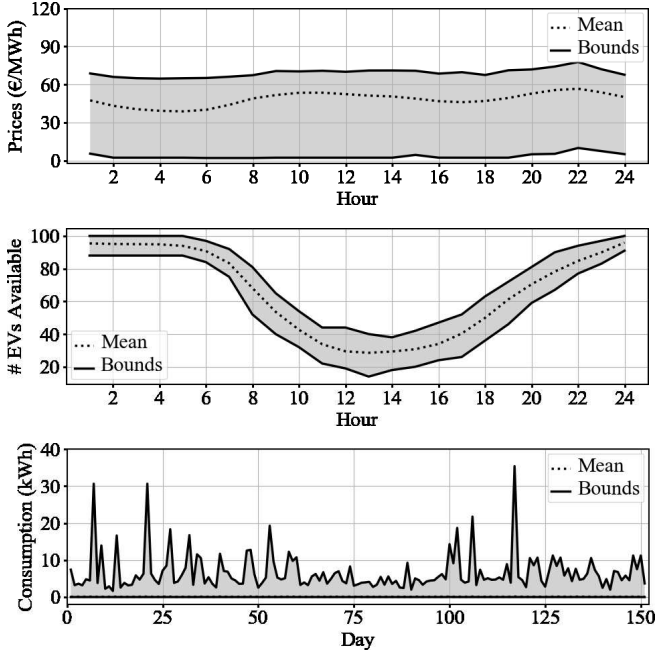}}
    \vspace{-0.4cm}
    \caption{Electricity prices, number of EVs available, and consumption due to transportation for a fleet of 100 EVs.}
    \label{figR_2}
\end{figure}

We assume a residential district aggregator with a fleet of 100 EVs. For the sake of simplicity, the technical parameters of each EV are identical. According to \cite{Technical_ZOE}, we consider that the maximum usable energy capacity of the EV batteries is 41.1 kWh, and thus, their minimum and maximum energy capacity are assumed to be 10 and 51.1 kWh, respectively; and that the energy rating per kilometer is 0.137 kWh/km. The maximum charging and discharging power is 7.4 kW; the round trip efficiency is 0.95; the battery cost is 70 $\textup{\euro}$/kWh; and the slope of the linear approximation of the battery life is -0.015625. Due to the lack of real-life data on the parameters associated with the driving patterns, we synthetically derive the availability profiles and energy required for transportation of EVs from the data collected by the National Household Travel Survey \cite{NHTS}, as described in \cite{8848991}. The electricity prices are provided by the ENTSO-e Transparency Platform \cite{ENTSOE} for the year 2018 in Spain. For the sake of illustration, Fig. \ref{figR_2} provides the input data related to historical records of electricity prices and travel patterns from January 1$^{st}$ till May 31$^{st}$, 2018. The upper plot in Fig. \ref{figR_2} shows the range of electricity prices, the plot in the middle represents the variability in the number of EVs available for a given hour, and the lower plot depicts the energy consumption due to motion for each day of the selected period. In addition, we consider that the penalty costs $C^p_1$ and $C^p_2$ are set to 2000 and 1000 $\textup{\euro}$/kWh, respectively. The capacity of the feeder is assumed unlimited, i.e. $P^G = 8000$ kW, unless stated otherwise. This is equivalent to disregarding the feeder capacity constraints in both the proposed and benchmark methodologies. Finally, daily simulations with hourly time steps have been run for four months spanning from February 1$^{st}$ till May 31$^{st}$.

Let us assume that the optimal market strategy is being determined for February 1$^{st}$ (Thursday). In order to focus on the impact of the uncertain availability of EVs, the electricity prices considered for all models are the average of the prices of the four previous days. The different models use the historical availability and consumption data as follows:
\begin{itemize}
    \item[-] The deterministic model is solved using the average availability and consumption of the previous four Thursdays.
    \item[-] The stochastic model is solved using four equiprobable scenarios containing the availability and consumption for each of the previous four Thursdays.
    \item[-] The set $\phi_{v}$ required by the proposed hierarchical model is computed using the observed availability for the previous four Thursdays, as explained in Subsection \ref{Hierarchical Formulation}. The  consumption of each EV is set to the average value for such days.
    \end{itemize}



The simulations have been run on a Linux-based server with one CPU clocking at 2.6 GHz and 2 GB of RAM using CPLEX 12.6.3 \cite{cplex} under Pyomo 5.2 \cite{pyomo}.

\subsection{Results for the Base Case}
\label{Results for the Base Case}

Table \ref{tab:results_basecase} summarizes the results from the proposed approach HF and the benchmarks DF and SF. We can find the following metrics for the whole optimization horizon: The total day-ahead cost $TC^{DA}$, the day-ahead purchase cost $C^{DA}$, the day-ahead battery degradation cost $D^{DA}$, the day-ahead sale revenue $R^{DA}$, the total energy bought and sold, i.e. $E^{B}$ and $E^{S}$, and deviations from both the energy balance $s_{FP}$ and the energy sale cleared in the day-ahead market $E^{-}_{FP}$. The latter metrics are obtained from the feasibility problem described in Section \ref{Comparison methodologies}. 

\begin{table}[tbp]
\caption{Results -- Base Case}
\vspace{-0.2cm}
\begin{center}
\begin{tabular}{ccccc}
 \hline
 \multicolumn{2}{c}{\textbf{Metric}} &    \textbf{DF}  &  \textbf{SF} & \textbf{HF}\\
\hline
\multicolumn{2}{l}{$TC^{DA}$ (\euro)}  &2282.4  &2708.4  &2888.4    \\
\multicolumn{2}{r}{$C^{DA}$ (\euro)}   &5875.2    &5709.6  &4314.0  \\
\multicolumn{2}{r}{$D^{DA}$ (\euro)}   &1686.0  &1548.0    &1143.6  \\
\multicolumn{2}{r}{$R^{DA}$ (\euro)}   &5278.8  &4549.2  &2569.2    \\
\hline
\multicolumn{2}{c}{$E^{B}$ (MWh)}   &162.2  &155.1 &114.2  \\
\multicolumn{2}{c}{$E^{S}$ (MWh)}   &96.5   &83.1  &47.7   \\
\hline
\multicolumn{2}{c}{$s_{FP}$ (MWh)}    &10.3     &4.7   &4.0  \\
\multicolumn{2}{c}{$E^{-}_{FP}$ (MWh)} &13.4    &1.2   &0.4 \\
\hline
\end{tabular}
\label{tab:results_basecase}
\end{center}
\end{table}

The total cost of the optimal day-ahead operation attained by the proposed approach HF amounts to \euro2888.4. Such a strategy leads to 4.0 MWh of energy deviations $s_{FP}$ for the EV-batteries, whereas the deviations $E^{-}_{FP}$ from the minimum value of the energy sold are equal to 0.4 MWh. The deterministic and stochastic methods (DF and SF) achieve 21.0\% and 6.2\% reduction in total day-ahead cost, in that order, compared to the proposed HF, as these models are less conservative than HF. However, the total day-ahead costs attained by DF and SF are reduced at the expense of (i) increasing the energy deviations from EV-batteries by 157.5\% and 17.5\%, respectively; and (ii) substantially increasing the deviations from the minimum value of the energy sold from 0.4 MWh for HF to 13.4 and 1.2 MWh for the DF and SF approaches, respectively. The high energy deviations attained by the deterministic method DF compared to HF are expected since DF is unable to capture the uncertainty in driving patterns, as opposed to SF and HF. Finally, when it comes to the degradation costs of the EV-batteries, DF and SF are more costly than the hierarchical model HF, since both DF and SF lead to more aggressive market strategies in terms of arbitrage compared to HF. As a consequence, both DF and SF buy more energy than HF, i.e. the EVs are scheduled to charge more power, thus the day-ahead purchase cost increases. In addition, the aggregator sells more energy when using both DF and SF (the EVs are scheduled to discharge more power) than HF, and thus leading to greater sale revenues. As a result, the degradation costs attained by DF and SF increase due to a rise in the number of charging/discharging cycles in the EV-batteries.

\begin{figure}[h]
\centerline{\includegraphics[scale=0.35]{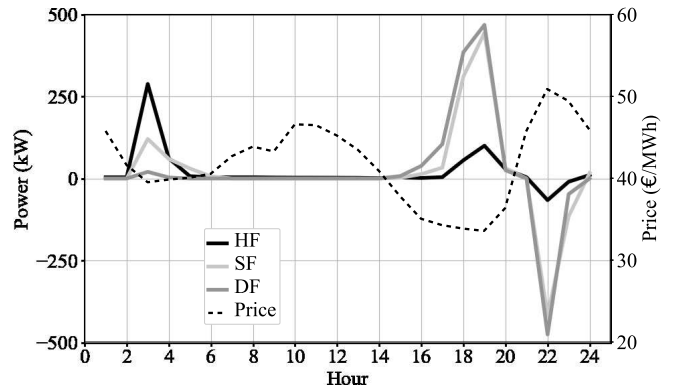}}
\vspace{-0.4cm}
\caption{Power bought and sold (left y-axis) and electricity prices (right y-axis) for each hour of day 89 for the base case.}
\label{figR_2_PP100}
\end{figure}

Fig. \ref{figR_2_PP100} provides the hourly power profiles, i.e. the power bought and sold in the day-ahead electricity market, for day 89 and each of the three methodologies. In this figure, we also show the electricity prices to better understand the periods when the aggregator purchases or sells power. As can be observed, although HF buys more energy at the beginning of the day than DF and SF, the proposed approach is much more conservative at the end of the day. Thus, HF schedules less power to be discharged into the grid. These power profiles are aligned with the aggregated results collated in Table \ref{tab:results_basecase}.

The daily average computational time for this base case with 100 EVs is 11.6s, 12.3s and 2.7s for HF, SF, and DF, in that order. 
    
\subsection{Impact of the Parameter $K_v$}
\label{Impact of the Parameter K_v}
The parameter $K_{v}$ sets the minimum time periods that the EV $v$ is estimated to be available. This way, the unavailability for charging and discharging of that EV is considered higher as $K_{v}$ decreases. In order to assess the impact of $K_{v}$, we have run two extra cases in which the value of $K_{v}$ is increased/decreased by five hours whenever feasible. Let us denote them as $K_v - 5$ and $K_v + 5$, respectively. Table \ref{tab:results_kv} summarizes the results for those cases along with the base case (simply detoned as $K_{v}$). When decreasing the value of $K_{v}$ up to 5 periods, the energy deviations $s_{FP}$ are reduced by 22.5\% with respect to the base case. This happens because the aggregator buys more energy, i.e. the purchase costs increase by 11.9\%, and sells less energy, i.e. the sale revenues decrease by 3.0\%, to hedge against the increased unavailability of the EV-fleet. As a consequence, the deviations $E^{-}_{FP}$ from the minimum value of power sold are equal to 0.0 MWh. As expected, this strategy leads to a more robustified solution for the aggregator. Conversely, when the EV has more time periods to be available, the aggregator behaves as a risk-taker with the goal of maximizing profits. Consequently, the energy deviations $s_{FP}$ increase by 20.0\% and $E^{-}_{FP}$ by 75\% from the base case.
\begin{table}[tbp]
\caption{Impact of the Parameter $K_v$}
\vspace{-0.2cm}
\begin{center}
\begin{tabular}{ccccc}
 \hline
 \multicolumn{2}{c}{\textbf{Metric}} &    \textbf{$K_v$}  &  \textbf{$K_v - 5$} & \textbf{$K_v + 5$}\\
\hline
\multicolumn{2}{l}{$TC^{DA}$ (\euro)} &2888.4   &3477.6 &2686.8    \\
\multicolumn{2}{r}{$C^{DA}$ (\euro)}   &4314    &4827.6 &4206    \\
\multicolumn{2}{r}{$D^{DA}$ (\euro)}   &1143.6  &1142.4 &1159.2    \\
\multicolumn{2}{r}{$R^{DA}$ (\euro)}   &2569.2  &2492.4 &2678.4    \\
\hline
\multicolumn{2}{c}{$E^{B}$ (MWh)}   &114.2  &125   &113.4 \\
\multicolumn{2}{c}{$E^{S}$ (MWh)}   &47.7   &46.2  &50  \\
\hline
\multicolumn{2}{c}{$s_{FP}$ (MWh)}     &4.0      &3.1  &4.8 \\
\multicolumn{2}{c}{$E^{-}_{FP}$ (MWh)}  &0.4      &0.0  &0.7    \\
\hline
\end{tabular}
\label{tab:results_kv}
\end{center}
\end{table}

\subsection{Scalability Results}
\label{Impact of the Feeder Capacity}
To test the scalability of the proposed methodology, we assume a fleet including 1000 EVs in the residential district, which we believe is large enough for this purpose. As similarly done for the base case, we provide the results for the proposed approach HF and the benchmarks DF and SF in Table \ref{tab:results_1000EVs} by considering an unconstrained feeder. From economic and technical standpoints, similar conclusions to the ones obtained for the base case can be drawn from Table \ref{tab:results_1000EVs} in the network-unconstrained case. The conservatism of HF leads to an increase in the total day-ahead cost of 3.5\% and 11.0\% compared to the results attained by SF and DF, however the energy deviations from EV-batteries decrease by 9.2\% and 63.6\%, respectively. In addition, the deviations from the minimum value of the energy sold decline to 0 when using HF. The robust decisions taken by HF are obviously reflected in the energy cleared from the day-ahead market. We can clearly see a reduction of \textit{circa} 25\% in the energy bought and almost 50\% in the energy sold for HF compared to the benchmarks.

From a computational standpoint, the results from Table \ref{tab:results_1000EVs} are achieved after 90.5s, 121.7s and 27.8s for HF, SF, and DF, in that order. The computational burden of DF is low compared to HF and SF because DF ignores the uncertainty in the driving patterns. Besides, HF is substantially faster than SF when increasing the EV-fleet size due to the complexity of the scenario-based modeling by SF. This demonstrates that HF scales well with the number of EVs in the fleet. Besides, the fact that our approach  encodes the uncertainty in the availability of each individual EV by way of a few intuitive parameters facilitates its implementation in real life, since those parameters could be directly provided by the EV users themselves through a home energy management system.

\begin{table}[h]
    \caption{Results -- Case 1000 EVs}
    \label{tab:results_1000EVs}
    \vspace{-0.2cm}
    \begin{center}
    \begin{tabular}{ccccc}
     \hline
     \multicolumn{2}{c}{\textbf{Metric}} &    \textbf{DF}  &  \textbf{SF} & \textbf{HF}\\
    \hline
    \multicolumn{2}{l}{$TC^{DA}$ (\euro)}  &23347.2  &25032.0  &25920.0    \\
    \multicolumn{2}{r}{$C^{DA}$ (\euro)}   &58857.6  &56074.8  &43486.8  \\
    \multicolumn{2}{r}{$D^{DA}$ (\euro)}   &16862.4  &15789.6  &11596.8  \\
    \multicolumn{2}{r}{$R^{DA}$ (\euro)}   &53372.8  &46832.4  &29163.6   \\
    \hline
    \multicolumn{2}{c}{$E^{B}$ (MWh)}   &1622.8  &1528.2    &1154.2  \\
    \multicolumn{2}{c}{$E^{S}$ (MWh)}   &957.4   &856.2     &481.3   \\
    \hline
    \multicolumn{2}{c}{$s_{FP}$ (MWh)}    &114.2     &45.8   &41.6  \\
    \multicolumn{2}{c}{$E^{-}_{FP}$ (MWh)} &129.5    &6.6   &0.0 \\
    \hline
    \end{tabular}
    \end{center}
    \end{table}

In this section, we also analyze the impact of the feeder capacity $P^G$ on the EV-fleet operation, which may simulate a more realistic setup. To emphasize such an impact, we assume the same fleet of 1000 EVs and reduce $P^G$ by 25\%, 50\% and 75\%, whose results are shown in Table \ref{tab:results_feeder}. First, $E^{-}_{FP}$ is equal to $0$ for all cases, i.e. there are only deviations from the energy balance of the vehicles' batteries regardless of $P^G$. As can be seen, the total day-ahead cost increases by 12.5, 13.2, and 17.5\% as the feeder capacity is reduced by 25, 50, and 75\%, respectively, compared to the unlimited case. Restricting the feeder capacity makes the aggregator to perform a less aggressive arbitrage strategy, i.e., the vehicles' charging power is more and more limited, thus leading to a decrease in its discharging power. This can be translated into lower sale revenues compared to the unlimited case. Besides, the energy deviations from EV-batteries slightly decrease up to 9.1\% for the most restrictive case. This reduction stems from the fact that the energy bought in the day-ahead market is redistributed into a greater number of periods when lowering the value of $P^G$.
\begin{table}[tbp]
\caption{Impact of the Feeder Capacity}
\vspace{-0.2cm}
\begin{center}
\begin{tabular}{cccccc}
 \hline
 \multicolumn{2}{c}{\textbf{Metric}} &    \textbf{0\%}  &  \textbf{25\%} & \textbf{50\%}  &  \textbf{75\%}\\
\hline
\multicolumn{2}{l}{$TC^{DA}$ (\euro)}  &25920    &29163.6  &29337.6  &30457.2\\
\multicolumn{2}{r}{$C^{DA}$ (\euro)}   &43486.8  &43482    &42512.4  &35048.4\\
\multicolumn{2}{r}{$D^{DA}$ (\euro)}   &11596.8  &11592    &11205.6  &8764.8\\
\multicolumn{2}{r}{$R^{DA}$ (\euro)}   &29163.6  &25910.4  &24380.4  &13356\\
\hline
\multicolumn{2}{c}{$E^{B}$ (MWh)}   &1154.2  &1153.9   &1120     &891.5\\
\multicolumn{2}{c}{$E^{S}$ (MWh)}   &481.3   &481      &450.4    &245.2\\
\hline
\multicolumn{2}{c}{$s_{FP}$ (MWh)}  &41.6      &41.6   &39.9    &37.8   \\
\hline
\end{tabular}
\label{tab:results_feeder}
\end{center}
\end{table}

\section{Conclusion}
\label{sec:conclusion}
This paper presents a computationally efficient robust optimization approach to model the profit-maximizer aggregator's decisions on the EV-fleet operation. A hierarchical program is used to represent the aggregator's operational model in the upper level, whereas a series of lower-level problems computes the vehicles' availability profiles that are worst case in terms of battery draining and energy exchange with the market. The original program is then reformulated as a mixed-integer linear model thanks to the unimodularity of the system matrices and results from duality theory. 

From the numerical results, we can conclude that the proposed model leads to a robustified and notably safer operation of the EVs in terms of deviations from the energy balance of their batteries, i.e. it can achieve reductions around 9--15\% and 60--64\% (depending on the fleet size) compared to stochastic and deterministic models, in that order. Naturally, those reductions come at the expense of increasing the total trading costs in the day-ahead market by 3--7\% and 11--27\% (depending also on the fleet size) compared to stochastic and deterministic models, respectively. Besides, the computational speed of the proposed model is up to 25\% faster than its stochastic counterpart for a real-sized electric-vehicle fleet. Finally, the feeder capacity can severely impact the operation of the fleet by diminishing its arbitrage opportunities. Specifically, the energy deviations from vehicles' batteries slightly decrease up to 9\% at the expense of increasing the day-ahead cost by 17.5\% for a 1000-vehicle fleet in a network-congested case.

\bibliographystyle{ieeetr}

\end{document}